\documentclass[12pt]{amsart}
\usepackage{fullpage,url,amssymb}
 \usepackage[all]{xy} 

\DeclareFontEncoding{OT2}{}{} 
\newcommand{\textcyr}[1]{%
 {\fontencoding{OT2}\fontfamily{wncyr}\fontseries{m}\fontshape{n}\selectfont #1}}
\newcommand{\Sha}{{\mbox{\textcyr{Sh}}}}

\usepackage{color}

\newcommand{\Cbar}{{\overline{C}}}

\newcommand{\Adual}{{A^\vee}}

\newcommand{\F}{{\mathbb F}}
\newcommand{\G}{{\mathbb G}}

\newcommand{\Q}{{\mathbb Q}}
\newcommand{\R}{{\mathbb R}}
\newcommand{\Z}{{\mathbb Z}}

\newcommand{\kbar}{{\overline{k}}}

\newcommand{\pp}{{\mathfrak p}}

\newcommand{\calA}{{\mathcal A}}

\newcommand{\calC}{{\mathcal C}}

\newcommand{\OO}{{\mathcal O}}

\DeclareMathOperator{\Gal}{Gal}

\DeclareMathOperator{\Br}{Br}

\DeclareMathOperator{\PIC}{\bf Pic}

\newcommand{\et}{{\operatorname{et}}}

\newcommand{\isom}{\simeq}

\newcommand{\intersect}{\cap} 
\newcommand{\tensor}{\otimes}

\newtheorem{theorem}{Theorem}[section]

\newtheorem{conjecture}[theorem]{Conjecture}

\theoremstyle{definition}

\newtheorem{question}[theorem]{Question}

\theoremstyle{remark}
\newtheorem{remark}[theorem]{Remark}

\usepackage[alphabetic,lite]{amsrefs} 

\begin{document}

\title[Brauer-Manin obstruction for curves]{Heuristics for the Brauer-Manin obstruction for curves}
\author{Bjorn Poonen}
\thanks{This research was supported by NSF grant DMS-0301280
          and a Packard Fellowship.  The author thanks the Institut
        Henri Poincar\'e and the Isaac Newton Institute for their hospitality.}
\address{Department of Mathematics, University of California, 
        Berkeley, CA 94720-3840, USA}
\email{poonen@math.berkeley.edu}
\urladdr{http://math.berkeley.edu/\~{}poonen}
\date{July 9, 2005}

\begin{abstract}
We conjecture that if $C$ is a curve of genus $>1$
over a number field $k$
such that $C(k)=\emptyset$,
then a method of Scharaschkin 
(equivalent to the Brauer-Manin obstruction in the context of curves)
supplies a proof that $C(k)=\emptyset$.
As evidence, we prove a corresponding statement
in which $C(\F_v)$ is replaced by a random subset of
the same size in $J(\F_v)$ for each residue field $\F_v$
at a place $v$ of good reduction for $C$,
and the orders of Jacobians over finite fields 
are assumed to be smooth (in the sense of having only small prime divisors)
as often as random integers of the same size.
If our conjecture holds, and if Shafarevich-Tate groups are finite,
then there exists an algorithm to decide whether a curve over $k$
has a $k$-point, and the Brauer-Manin obstruction to the Hasse principle
for curves over the number fields is the only one.
\end{abstract}

\maketitle

\section{Setup}\label{S:setup}
Let $k$ be a number field.
Fix an algebraic closure $\kbar$ of $k$, and let $G=\Gal(\kbar/k)$.
Let $C$ be a curve of genus $g$ over $k$.
(In this paper, curves are assumed to be 
smooth, projective, and geometrically integral.)
Let $\Cbar = C \times_k \kbar$.
Let $J$ be the Jacobian of $C$, 
which is an abelian variety of dimension $g$ over $k$.
Assume that $\Cbar$ has a $G$-invariant line bundle of degree $1$:
this gives rise to a $k$-morphism $C \to J$,
and it is an embedding if $g>0$.
Let $S_C$ be the set of finite primes $v$ of good reduction for $C$.
Similarly define $S_A$ for any abelian variety $A$ over $k$.
We have $S_C \subseteq S_J$.

\section{Determining the set of rational points}
Suppose that generators of the Mordell-Weil group $J(k)$ are known.
Then $C(k)$ equals the set of points in $J(k)$ 
that lie on the subvariety $C$.
We would like to know whether $C(k)$ can be calculated,
especially in the case $g>1$ in which $C(k)$ is guaranteed to be finite 
by~\cite{Faltings1983}.

If $J(k)$ is finite, 
then in principle we can list its elements 
and check which of them lie on $C$.
On the other hand, if $J(k)$ is infinite, 
it can be very difficult to decide which points of $J(k)$ lie on $C$.

\section{Chabauty's approach}
One approach, due to C.~Chabauty~\cite{Chabauty1941},
works in $J(k_v)$, where $k_v$ is the completion of $k$
at a nonarchimedean place $v$.
Chabauty observed that the closure $\overline{J(k)}$ of $J(k)$ in $J(k_v)$
is an analytic submanifold of $J(k_v)$,
and proved that {\em if the rank of $J(k)$ is less than $g$},
the subset of points in $\overline{J(k)}$ lying on $C$ is finite;
in this case, as explained by R.~Coleman~\cite{Coleman1985chabauty}, 
one gets an effective upper bound for $\#C(k)$.
Often one can even determine $C(k)$ explicitly.
But the dimension hypothesis is not always satisfied,
and even when it is, the upper bound on $\#C(k)$ may fail to be sharp.

\section{Scharaschkin's approach}
A more recent approach, 
suggested by V.~Scharaschkin~\cite{Scharaschkin2004preprint},
tries to find which points of $J(k)$ lie on $C$ modulo $\pp$ 
for many primes $\pp$.
More precisely, he proposes the following method
for proving that $C(k)$ is empty.
Choose a finite subset $S \subset S_C$.
Let $\F_v$ be the residue field at $v$.
Then we have a commutative square
\[
\xymatrix{
C(k) \ar@{..>}[r] \ar@{..>}[d] & \prod_{v \in S} C(\F_v) \ar[d] \\
J(k) \ar[r] & \prod_{v \in S} J(\F_v). \\
}
\]
If the images of the solid arrows in $\prod_{v \in S} J(\F_v)$
do not intersect, then $C(k)$ is empty.

\begin{remark}\label{R:B-M obstruction}
For this paragraph we assume that the 
Shafarevich-Tate group $\Sha(J)$ is finite,
or at least that its maximal divisible subgroup is trivial.
Scharaschkin~\cite{Scharaschkin2004preprint} proved that then 
the potential obstruction to the existence of $k$-points
described above
is part of the Brauer-Manin obstruction.
More precisely, 
he showed that if we consider the product of $J(k_v)$ 
(modulo its connected component if $v$ is archimedean) 
over all places $v$ 
instead of a product of only $J(\F_v)$ over only places of good reduction,
then we recover exactly the Brauer-Manin obstruction.

For the connection of the Brauer-Manin obstruction
to the information on rational points
obtained from finite \'etale covers, see~\cite{Stoll2005preprint}.
\end{remark}

\section{A conjecture and its implications}

We conjecture the following, based on heuristics to be explained later.

\begin{conjecture}\label{C:do not intersect}
Let $C \to J$ be as in Section~\ref{S:setup}, with $g>1$.
If $C(k)=\emptyset$,
then there exists a finite subset $S \subset S_C$
such that the images of $J(k)$ and $\prod_{v \in S} C(\F_v)$ 
in $\prod_{v \in S} J(\F_v)$
do not intersect.
\end{conjecture}

The importance of Conjecture~\ref{C:do not intersect} is given
by the following result:

\begin{theorem}\label{T:important}
Assume Conjecture~\ref{C:do not intersect}.
Assume also that Shafarevich-Tate groups of Jacobians of
curves over number fields are finite.
Then
\begin{enumerate}
\item[(a)]
There is an algorithm that takes as input a number field $k$ and a curve $C$ over $k$, and decides whether $C$ has $k$-point.
\item[(b)]
The Brauer-Manin obstruction to the Hasse principle
is the only obstruction to the existence of a $k$-point
on a curve $C$ over a number field $k$.
\end{enumerate}
\end{theorem}

\begin{proof}
For details on how elements of $k$ and curves over $k$
can be represented, see~\cite{Baker-et-al2004preprint}.

Before proceeding, we recall a few well-known facts:
\begin{enumerate}
\item[(i)]
There exists an algorithm
for deciding whether a smooth projective variety $X$ over $k$
has a point over every completion $k_v$.
(Sketch of proof: For all nonarchimedean primes $v$ of sufficiently
large norm, the Weil conjectures and Hensel's lemma imply
that $X$ automatically has a $k_v$-point.
One can test the remaining $v$ individually, again using Hensel's lemma
in the nonarchimedean case.)
\item[(ii)]
If $X$ is a torsor of an abelian variety $A$ over a number field $k$,
and if $\Sha(A)$ is finite,
then there exists an algorithm to decide whether $X$ has a $k$-point,
and the Brauer-Manin obstruction to the Hasse principle is the
only one for $X$.
(Sketch of proof: 
By the previous fact, we may assume $X(k_v) \ne \emptyset$ for all $v$,
so $X$ represents an element of $\Sha(A)$.
If we search for $k$-points on $X$ by day,
and perform higher and higher descents on $A$ by night,
we will eventually decide whether $X$ has a $k$-point,
assuming the finiteness of $\Sha(A)$,
even if we are not given a bound on $\#\Sha(A)$.
It remains to show that if $X$ has no $k$-point,
then there is a Brauer-Manin obstruction.
Under our assumption that $\Sha(A)$ is finite,
the Cassels-Tate pairing
\[
        \langle \;,\; \rangle\colon \Sha(A) \times \Sha(\Adual) \to \Q/\Z
\]
is nondegenerate.
If $X$ has no $k$-point, then $X$ corresponds to a nonzero element
of $\Sha(A)$, so there is a torsor $Y$ of the dual abelian variety $\Adual$
such that $\langle X,Y \rangle \ne 0$.
By \cite{Manin1971}*{Theorem~6},
there is an element $y \in \Br X$ related to $Y$
such that the Brauer pairing of $y$ with every adelic point on $X$
gives the value $\langle X,Y \rangle \ne 0$,
so $X$ has a Brauer-Manin obstruction.)
\end{enumerate}

We now return to our problem.
By \cite{Baker-et-al2004preprint}*{Lemma~5.1(1)}, we can compute
the genus $g$ of $C$, so we may break into cases according to
the value of $g$.

If $g=0$, then $C$ satisfies the Hasse principle.
Thus to test for the existence of a $k$-point,
it suffices to use fact~(i) above.

If $g=1$, then $C$ is a torsor of its Jacobian $J$,
so it suffices to use~(ii).

{}From now on, we suppose $g \ge 2$.
By~(i),
we may reduce to the case that $C(k_v)$ is nonempty for every $v$.
Let $X:=\PIC^1_{C/k}$ be the variety 
parametrizing degree-$1$ line bundles on $C$.
Thus $X$ is a torsor of the Jacobian $J$.
We have a canonical injection from $C$ to $X$
taking a point $c \in C$ to the class of the associated degree-$1$ divisor.
By~(ii),
we can check whether $X$ has a $k$-point;
if not, then $C$ has no $k$-point,
and there is a Brauer-Manin obstruction for $X$,
which pulls back to a Brauer-Manin obstruction for $C$.
Thus from now on, we may assume that $X$ has a $k$-point.
In other words, $\Cbar$ has a $G$-invariant line bundle of degree $1$.
Such a $G$-invariant line bundle can be found by a search,
and it allows us to identify $X$ with $J$.
We now can search for $k$-points on $C$ each day,
while running Scharaschkin's method using the first $r$ primes in $S_C$
for larger and larger $r$ each night.
Conjecture~\ref{C:do not intersect} implies that one of these two processes
will terminate.
Thus there exists an algorithm for deciding 
whether $C$ has a $k$-point.
Moreover, as mentioned in Remark~\ref{R:B-M obstruction},
assuming finiteness of $\Sha(J)$,
if Scharaschkin's method proves the nonexistence of $k$-points,
then there is a Brauer-Manin obstruction.
\end{proof}

\begin{remark}
If one knows that the Brauer-Manin obstruction to the Hasse principle
is the only one for a smooth projective variety $X$,
that in itself lets one determine whether $X$ has a $k$-point,
in principle, as we will explain in the following paragraph.
This gives an alternative approach to Theorem~\ref{T:important}(a),
based on part~(b).

By an unpublished result of O.~Gabber, reproved by A.~J.~de~Jong,
each element of the cohomological Brauer group $\Br X:=H^2_{\et}(X,\G_m)$ 
can be represented by an Azumaya algebra $\calA$, 
i.e., a locally free $\OO_X$-algebra that is \'etale locally
isomorphic to a finite-dimensional matrix algebra.
Each $\calA$ can be described by a finite amount of data:
\begin{enumerate}
\item
a covering of $X$ by finitely many Zariski open subsets $X_i$ 
such that $\calA|_{X_i}$ is free as an $\OO_{X_i}$-module,
\item
the multiplication table for $\calA|_{X_i}$ 
with respect to a chosen $\OO_{X_i}$-basis, for each $i$,
\item
the change-of-basis map 
on the intersection $X_i \intersect X_j$, for each $i$ and $j$,
\item
a covering $U_i \to X_i$ in the \'etale topology, for each $i$,
\item
a positive integer $r_i$ and an $\OO_{U_i}$-algebra isomorphism
$\calA \tensor_{\OO_X} \OO_{U_i} \isom M_{r_i}(\OO_{U_i})$, 
for each $i$.
\end{enumerate}
Moreover, given such data, one can easily check whether it actually defines
an Azumaya algebra.
Therefore there is a (very inefficient) algorithm 
which when left running forever,
eventually produces all Azumaya algebras over $X$,
each possibly occurring more than once,
simply by enumerating and checking each possible set of data.
Now, we search for $k$-points by day and generate Azumaya algebras by night, 
calculating at the end of each night
whether the Azumaya algebras generated so far give an obstruction.
\end{remark}

\begin{remark}
It is not clear whether Conjecture~\ref{C:do not intersect}
and the finiteness of Shafarevich-Tate groups
imply the existence of an algorithm for {\em listing}
all $k$-points on a given curve $C$ of genus $g \ge 2$ over $k$.
For the listing problem, applying Chabauty's method to finite \'etale
covers seems more promising: see~\cite{Stoll2005preprint} for an analysis
of the situation.
\end{remark}

\section{Computational evidence for the conjecture}

E.~V.~Flynn~\cite{Flynn2004} 
has developed an implementation of Scharaschkin's method
for genus-$2$ curves over $\Q$.
He tested $145$ such curves defined by equations with small integer
coefficients, having $\Q_p$-points for all $p \le \infty$,
but having no $\Q$-point with $x$-coordinate 
of height less than $10^{30}$.
These were grouped according to the rank $r$ of the Jacobian.
In all cases with $r \le 1$, he successfully showed that
there was a Brauer-Manin obstruction.
In most cases with $r=2$, a Brauer-Manin obstruction was found,
and all the unresolved $r=2$ cases were later resolved by an improved 
implementation of M.~Stoll.
The remaining cases had $r=3$, and a few of these were resolved;
it was unclear from the computation whether the remaining ones
could be resolved by a longer computation:
the combinatorics quickly became prohibitive.

\section{Theoretical evidence for the conjecture}

Here we give a heuristic analysis of Conjecture~\ref{C:do not intersect}.

Recall that if $B \in \R_{>0}$, 
an integer is called {\em $B$-smooth} if all its prime factors are $\le B$.
For any fixed $u \in (0,1)$, 
the fraction of integers in $[1,B]$ that are $B^u$-smooth
tends to a positive constant as $B \to \infty$ \cite{DeBruijn1951}.

As our main evidence for Conjecture~\ref{C:do not intersect},
we prove a modified version of it in which $C(\F_v)$ is modelled
by a random subset of $J(\F_v)$ of the same order,
and in which we assume that the integer $\# J(\F_v)$ is 
as smooth as often as a typical integer of its size.
The smoothness assumption is formalized in the following:

\begin{conjecture}\label{C:smoothness}
Let $A$ be an abelian variety over a number field $k$,
and let $u \in (0,1)$.
Then
\[
\limsup_{B \to \infty} \frac{\{v \in S_A : \#\F_v \le B \text{ and $\#A(\F_v)$ is $B^u$-smooth }\}}
{\{v \in S_A : \#\F_v \le B \}} > 0.
\]
\end{conjecture}

Let $g=\dim A$.  
If $\#A(\F_v)$ behaves like a typical integer of its size,
which is about $(\#\F_v)^g \le B^g$,
then it should have a positive probability of being $B^u$-smooth,
since $B^u$ is a constant power of $B^g$.
If anything, $\#A(\F_v)$ can be expected to factor more than
typical integers its size, because of splitting of $A$ up to isogeny,
or because of biases in the probability of being divisible by small primes.
Thus Conjecture~\ref{C:smoothness} is reasonable.

We are now ready to state our main result giving evidence 
for Conjecture~\ref{C:do not intersect}.

\begin{theorem}
Let $C \to J$ be as in Section~\ref{S:setup}, with $g>1$.
Assume Conjecture~\ref{C:smoothness} for $J$.
For each prime $v \in S_C$,
let $\calC(\F_v)$ be a random subset of $J(\F_v)$ of size $\#C(\F_v)$;
we assume that the choices for different $v$ are independent.
Then with probability $1$, there exists a finite subset $S \subseteq S_C$
such that the images of $J(k)$ and $\prod_{v \in S} \calC(\F_v)$
in $\prod_{v \in S} J(\F_v)$
do not intersect.
\end{theorem}

\begin{proof}
It suffices to find $S$ such that the probability that the images intersect
is arbitrarily small.

Given $B>0$, let $S=S(B)$ be the set of $v \in S_C$
such that $\# \F_v \le B^2$ and $\#J(\F_v)$ is $B$-smooth.
Because we have assumed Conjecture~\ref{C:smoothness} for $J$,
there exists $c>0$ such that for arbitrarily large
$B \in \R_{>0}$ (the square root of the $B$ occurring
in Conjecture~\ref{C:smoothness}), 
the set $S$ contains at least a fraction $c$ of the primes $v \in S_C$
with $\#\F_v \le B^2$.

Let $\pi(x)$ be the number of rational primes $\le x$.
The prime number theorem says that 
$\pi(x) = (1+o(1)) x/\log x$ as $x \to \infty$.
For $v \in S$, the Weil conjectures give 
$\#J(\F_v) \le O\left((\#\F_v)^{2g}\right) \le B^{2g+o(1)}$ as $B \to \infty$
so by $B$-smoothness, 
the least common multiple $L$ of $\# J(\F_v)$ for $v \in S$ satisfies
\begin{align*}
L &\le \prod_{\text{primes $p \le B$}} p^{\lfloor \log_p B^{2g+o(1)} \rfloor} \\
        &\le \prod_{\text{primes $p \le B$}} B^{2g+o(1)} \\
\log L &\le \pi(B) (2g+o(1)) \log B \\
        &=  (2g+o(1)) B.
\end{align*}
Suppose that $J(k)$ is generated by $r$ elements.
Every element of $\prod_{v \in S} J(\F_v)$ has order dividing $L$,
so the order of the image $I$ of $J(k)$ in $\prod_{v \in S} J(\F_v)$ is at most
\[
        L^r \le \exp\left((2g+o(1)) r B \right).
\]

The probability that a fixed element of a set of size $n$ belongs
to a random subset of size $m$ is $m/n$,
so the probability $P$ that a fixed element of $I$ belongs to the image $I'$
of $\prod_{v \in S} \calC(\F_v) \to \prod_{v \in S} J(\F_v)$
satisfies
\begin{align*}
P       &= \prod_{v \in S} \frac{\#C(\F_v)}{\#J(\F_v)} \\
        &\le \prod_{v \in S} \frac{(\#\F_v)^{1+o(1)}}{(\#\F_v)^{g+o(1)}} \\
\log P  &= (1-g+o(1)) \sum_{v \in S} \log \#\F_v,
\end{align*}
where again we have used the Weil conjectures.
The number of primes $v$ of $k$ with $\#\F_v \le B^2$
is at least the number of rational primes $p \le B^2$
that split completely in $k$, which is asymptotically 
$\ge c_1 \pi(B^2)$.
(We use $c_1,c_2,\dots$ to denote positive constants independent of $B$.)
Since $S$ contains a positive fraction of these $v$,
and since the numbers $\#\F_v$ are powers of primes, with
each prime occurring at most $[k:\Q]$ times, 
we find that $\sum_{v \in S} \log \#\F_v$ is at least
the sum of $\log p$ for the first $c_2 \pi(B^2)$ primes.
By the prime number theorem, the $n$-th prime is $(1+o(1)) n \log n$,
so 
\begin{align*}
        \sum_{v \in S} \log \#\F_v 
        &\ge \sum_{n=1}^{c_2 \pi(B^2)} \log((1+o(1)) n \log n) \\
        &\ge \sum_{n=c_2 \pi(B^2)/2}^{c_2 \pi(B^2)} \log((1+o(1)) n \log n) \\
        &\ge c_3 \pi(B^2) \log \pi(B^2) \\
        &\ge c_4 B^2
\end{align*}
as $B \to \infty$.
Since $1-g<0$, we get
\begin{align*}
        \log P &\le -c_5 B^2 \\
        P &\le \exp(-c_5 B^2).
\end{align*}
Thus the probability that $I$ intersects $I'$ is at most
\[
        \#I \cdot P \le \exp\left((2g+o(1)) r B \right) \cdot \exp(-c_5 B^2),
\]
which tends to $0$ as $B \to \infty$, as desired.
\end{proof}

\begin{remark}
For the case where $C(k)$ is nonempty,
we would have liked to analyze a refined heuristic
that reflects the existence of the $k$-points.
Namely, suppose that $\calC(\F_v)$ is a random subset of $J(\F_v)$
of size $\#C(\F_v)$ chosen 
{\em subject to the constraint that it 
contains the image of $C(k)$ in $J(\F_v)$}.
We then expect that with probability~$1$,
the only points of $J(k)$ whose image in $\prod_{v \in S_C} J(\F_v)$
lie in the image of $\prod_{v \in S_C} \calC(\F_v)$ are those in $C(k)$.
But we were unable to prove this,
even assuming Conjecture~\ref{C:smoothness}.
\end{remark}

\begin{question}
Is it true more generally that if $X$ is a closed subvariety 
of an abelian variety $A$ over a number field $k$,
and $S$ is a density-$1$ set of primes of good reduction for $X$ and $A$,
then the intersection of the closure of the image of $A(k)$ in
$\prod_{v \in S} A(\F_v)$
with $\prod_{v \in S} X(\F_v)$
equals the closure of the image of $X(k)$ in $\prod_{v \in S} A(\F_v)$?
One could also ask the question with $\prod_{v \in S} A(\F_v)$
replaced by the product of $A(k_v)$ (modulo its connected component) 
over all $v$.
We expect a positive answer;
this together with finiteness of $\Sha(A)$
would imply that the Brauer-Manin
obstruction to the Hasse principle is the only one 
for such~$X$: cf.~Remark~\ref{R:B-M obstruction}.
\end{question}

\section*{Acknowledgements} 

I thank Jean-Louis Colliot-Th\'el\`ene and Michael Stoll
for several discussions about the Brauer-Manin obstruction 
during the Fall 2004 trimester at the Institut Henri Poincar\'e.

\end{document}